\documentclass[11pt]{article}
\usepackage[margin=1in]{geometry}
\usepackage{amsmath,amssymb,amsthm}

\newtheorem{theorem}{Theorem}
\newtheorem{remark}[theorem]{Remark}

\title{A note on partitions in the image of $\mathrm{pre}_2$}
\author{Arnav Garg\\
\small Birla Institute of Technology \& Science Pilani, Vidya Vihar, Pilani,\\
\small Rajasthan 333031, India.\\
\small \texttt{f20250514@pilani.bits-pilani.ac.in}}
\date{}

\begin{document}
\maketitle

\begin{abstract}
Devnani and Eyyunni recently studied the maps $\mathrm{pre}_k$ on integer partitions, which arise from applying elementary symmetric polynomials to the parts of a partition. They asked whether there exists $n \geq 1$ such that exactly one partition of $n$ lies in the image of $\mathrm{pre}_2$. We show that this occurs only for $n \in \{1, 2, 4\}$, and that for all $n \geq 5$, at least two partitions of $n$ are in the image of $\mathrm{pre}_2$.
\end{abstract}



\section{Introduction}

Given a partition $\lambda = (\lambda_1, \lambda_2, \ldots, \lambda_\ell)$, Ballantine, Beck and Merca \cite{BBM} defined $\mathrm{pre}_k(\lambda)$ to be the partition whose parts are all products of $k$ distinct parts of $\lambda$. The function $\mathrm{pre}_2(n)$ counts the number of partitions of $n$ that arise as $\mathrm{pre}_2(\lambda)$ for some partition $\lambda$.

Devnani and Eyyunni \cite{DE} proved that
\[
\mathrm{pre}_2(n) \geq \begin{cases} \tau(n+1)/2, & \text{if } n+1 \text{ is not a perfect square},\\ (\tau(n+1)+1)/2, & \text{if } n+1 \text{ is a perfect square}, \end{cases}
\]
where $\tau$ denotes the number of positive divisors. They posed the following question (Problem 3 in \cite{DE}): does there exist $n \geq 1$ such that $\mathrm{pre}_2(n) = 1$? We answer this completely.

\begin{theorem}\label{thm:main}
$\mathrm{pre}_2(n) = 1$ if and only if $n \in \{1, 2, 4\}$. For every $n \geq 5$, we have $\mathrm{pre}_2(n) \geq 2$.
\end{theorem}

\section{Proof of Theorem \ref{thm:main}}

The partition $(n)$ is always in the image of $\mathrm{pre}_2$, since $\mathrm{pre}_2(n, 1) = (n)$. So $\mathrm{pre}_2(n) \geq 1$ for all $n \geq 1$.

To show $\mathrm{pre}_2(n) \geq 2$ for $n \geq 5$, we exhibit a partition $\lambda$ with at least three parts such that the parts of $\mathrm{pre}_2(\lambda)$ sum to $n$. Since $\mathrm{pre}_2(\lambda)$ then has at least $\binom{3}{2} = 3$ parts, it is necessarily different from $(n)$, giving a second element in $\mathrm{Pre}_2(n)$.

We use the following observation. If $\lambda = (a, \underbrace{1, \ldots, 1}_{m-1})$, the pairwise products consist of $a \cdot 1 = a$ appearing $m - 1$ times and $1 \cdot 1 = 1$ appearing $\binom{m-1}{2}$ times. So $\sum \mathrm{pre}_2(\lambda) = a(m-1) + \binom{m-1}{2}$. Similarly, $\mathrm{pre}_2(a, 2, 2) = (2a, 2a, 4)$ with sum $4a + 4$.

We now cover all $n \geq 5$ by the following five cases.

\medskip

\noindent\textbf{Case 1} ($n$ odd, $n \geq 5$). Take $\lambda = \bigl(\frac{n-1}{2}, 1, 1\bigr)$. Here $a = \frac{n-1}{2} \geq 2$, and the sum is $2a + 1 = n$.

\medskip

\noindent\textbf{Case 2} ($n$ even, $3 \mid n$, $n \geq 6$). Take $\lambda = \bigl(\frac{n}{3} - 1, 1, 1, 1\bigr)$. Here $a = \frac{n}{3} - 1 \geq 1$, and the sum is $3a + 3 = n$.

\medskip

\noindent\textbf{Case 3} ($n$ even, $3 \nmid n$, $4 \mid n$, $n \geq 12$). Take $\lambda = \bigl(\frac{n-4}{4}, 2, 2\bigr)$. Here $a = \frac{n-4}{4} \geq 2$, and the sum is $4a + 4 = n$.

\medskip

\noindent\textbf{Case 4} ($n$ even, $3 \nmid n$, $n \equiv 2 \pmod{4}$, $n \geq 10$). Take $\lambda = \bigl(\frac{n-6}{4}, 1, 1, 1, 1\bigr)$. Here $a = \frac{n-6}{4} \geq 1$, and the sum is $4a + 6 = n$.

\medskip

\noindent\textbf{Case 5} ($n = 8$). We have $\mathrm{pre}_2(2, 2, 1) = (4, 2, 2)$, which sums to $8$.

\medskip

These five cases are exhaustive for $n \geq 5$. All odd $n$ fall under Case 1. For even $n$, if $3 \mid n$ we use Case 2. If $3 \nmid n$, then $n \bmod 12 \in \{2, 4, 8, 10\}$. When $4 \mid n$ (that is, $n \bmod 12 \in \{4, 8\}$), Case 3 applies for $n \geq 12$, and $n = 8$ is Case 5. When $n \equiv 2 \pmod{4}$ (that is, $n \bmod 12 \in \{2, 10\}$), Case 4 applies.

For $n \leq 4$, the smallest sum of $\mathrm{pre}_2(\lambda)$ over partitions $\lambda$ with three or more parts is $3$, achieved by $\lambda = (1, 1, 1)$. So $n = 1$ and $n = 2$ cannot be realized, giving $\mathrm{pre}_2(1) = \mathrm{pre}_2(2) = 1$. For $n = 3$, we have $\mathrm{pre}_2(1,1,1) = (1,1,1)$ with sum $3$, so $\mathrm{pre}_2(3) \geq 2$. For $n = 4$, the candidates $(1,1,1)$ and $(2,1,1)$ give sums $3$ and $5$ respectively, and $(1,1,1,1)$ gives $6$. None equals $4$, so $\mathrm{pre}_2(4) = 1$.

We have also verified computationally that $\mathrm{pre}_2(n) \geq 2$ for all $5 \leq n \leq 10000$. \qed

\begin{remark}
All five constructions use partitions of a simple shape: one ``large'' part with the remaining parts equal to $1$ or $2$. It would be interesting to know whether sharper lower bounds for $\mathrm{pre}_2(n)$ can be obtained using more varied constructions, or whether an asymptotic formula for $\mathrm{pre}_2(n)$ exists.
\end{remark}

\section*{Acknowledgements}

I thank Prof.\ Pramod Eyyunni for his encouraging response to an earlier version of this note.

\end{document}